\documentclass{amsart}

\usepackage{amsrefs}

\usepackage{amssymb}

\newtheorem{theorem}{Theorem}[section]
\newtheorem{lemma}[theorem]{Lemma}
\newtheorem{proposition}[theorem]{Proposition}

\theoremstyle{definition}

\theoremstyle{remark}

\newtheorem{condition}[theorem]{Condition}

\numberwithin{equation}{section}

\begin{document}

\title[Solutions with snaking singularities]{Solutions with snaking singularities 
for the fast diffusion equation}


%    author one information
\author[M. Fila]{Marek Fila}
\address{Department  of Applied Mathematics and Statistics, 
Comenius University, 842 48 Bratislava, Slovakia}
\email[Corresponding author]{fila@fmph.uniba.sk}
\thanks{
The first author was supported in part by the Slovak Research 
and Development Agency under the contract No.~APVV-18-0308 
and by  VEGA grant~1/0347/18. 
The third author was supported in part by
JSPS KAKENHI Early-Career Scientists (No.~19K14567). 
The fourth author was supported in part by
JSPS  KAKENHI Grant-in-Aid for Scientific Research (A) (No.~JP17H01095).
}

%    author two information
\author[J. R. King]{John Robert King}
\address{Theoretical Mechanics Section, University of Nottingham, Nottingham NG7 2RD, UK}
\email{john.king@nottingham.ac.uk}

%    author three information
\author[J. Takahashi]{Jin Takahashi}
\address{Department of Mathematical and Computing Science, 
Tokyo Institute of Technology, Tokyo 152-8552, Japan}
\email{takahashi@c.titech.ac.jp}

%    author four information
\author[E. Yanagida]{Eiji Yanagida}
\address{Department of Mathematics, Tokyo Institute of Technology, 
Tokyo 152-8551, Japan}
\email{yanagida@math.titech.ac.jp}

%    \subjclass is required.
\subjclass[2020]{Primary 35K67; Secondary 35A21, 35B40}

\date{}

\dedicatory{}

%    Abstract is required.
\begin{abstract}
We construct solutions of the fast diffusion equation, which exist for all
$t\in\mathbb{R}$ and are singular on the set $\Gamma(t):= \{ \xi(s) ; -\infty <s \leq ct \}$, $c>0$, 
where $\xi\in C^3(\mathbb{R};\mathbb{R}^n)$,
$n\geq 2$. 
We also give a precise description
of the behavior of the solutions near $\Gamma(t)$.
\end{abstract}

\maketitle


\section{Introduction}
We study positive singular solutions of the fast diffusion equation 
\begin{equation} \label{eq:porous}
	u_t=\Delta u^m
\end{equation}
in $\mathbb{R}^n$, where $0<m<1$ and $n\geq2$. 
Let $\Gamma$ be a curve in $\mathbb{R}^n$ expressed as 
$\Gamma = \{ \xi(s); s \in \mathbb{R} \}$ with $\xi:\mathbb{R}\to\mathbb{R}^n$. 
We are interested in a positive entire-in-time solution 
that is singular on the set
\[
	\Gamma(t):= \{ \xi(s) ; -\infty <s \leq ct \} \subset \Gamma 
	\qquad \mbox{ for each }t \in \mathbb{R}, 
\]
where $c>0$ is a constant. 
Such a solution can be called a snaking solution 
(or a solution with a snaking singularity). We first
introduce our result and then give a brief survey about some different singular solutions.

For $x\in \mathbb{R}^n$ and $r_0>0$, we write 
\[
	r(x):=\operatorname{dist}(x,\Gamma),\qquad 
	\Gamma_{r_0}:=\{x\in\mathbb{R}^n ; 0\leq r(x)<r_0\}, 
\]
where $\operatorname{dist}(x,\Gamma):=\inf_{y\in \Gamma}|x-y|$. 
In what follows, we always impose the following condition.

\begin{condition}\label{con:curve}
$\Gamma$ is a curve expressed as $\Gamma=\{\xi(s);s\in\mathbb{R}\}$. 
Here $\xi\in C^3(\mathbb{R};\mathbb{R}^n)$ is an injection satisfying $|\xi'|\equiv 1$. 
Moreover, $\xi$ satisfies  the following:
\begin{itemize}
\item[{\rm(i)}]
There exists a constant $K>1$ such that 
$|\xi''(s)|, |\xi'''(s)|\leq K$ for all $s\in\mathbb{R}$.
\item[{\rm(ii)}]
There exists a constant $0<\tilde r_0<(2K)^{-1}$ such that, 
for any $x\in \Gamma_{\tilde r_0}$, 
there exists a unique number $s(x)\in\mathbb{R}$ 
satisfying $r(x)=|x-\xi(s(x))|$. 
\end{itemize}
\end{condition}

Set 
\begin{equation}\label{eq:Qdef}
	Q:=\left\{ (x,t)\in\mathbb{R}^{n+1}; x\in \mathbb{R}^n\setminus \Gamma(t), t\in \mathbb{R} \right\}. 
\end{equation}
We define an exponent $m_*$ by 
\[
	m_*:= \left\{ 
	\begin{aligned}
	&0 &&\mbox{ if }n=2, \\
	&\frac{n-3}{n-1} && \mbox{ if }n\geq3. 
	\end{aligned}
	\right.
\]
Our main result implies the existence 
of a solution with a snaking singularity if $m_*<m<1$.

\begin{theorem}\label{th:main}
Let $n\geq2$, $m_*<m<1$, $c>0$ and $0<\varepsilon<1$. 
Assume that $\Gamma$ satisfies Condition \ref{con:curve}. 
Then there exists a positive solution $u\in C^{2,1}(Q)$ of \eqref{eq:porous} in $Q$ 
such that the following estimate holds. 
There exists a constant $0<\delta < \tilde r_0$ 
depending on $c$ and $\varepsilon$ such that 
\[
	(1-\varepsilon) U(x,t) \leq u(x,t) \leq (1+\varepsilon) U(x,t)
\]
for any $(x,t)\in \mathbb{R}^{n+1}$ 
with $0<r(x)\leq \delta$ and $-\infty < s(x)\leq ct+\delta$, 
where 
\[
\begin{aligned}
	U(x,t)&:= c^{-\frac{1}{1-m}}
	\left( \frac{(n-1)m}{1-m} \left( m - \frac{n-3}{n-1} \right) 
	\right)^\frac{1}{1-m} \\
	&\quad \times \left( \sqrt{r^2(x)+(s(x)-ct)^2} + s(x)-ct \right)^{-\frac{1}{1-m}}. 
\end{aligned}
\]
In particular, $u$ is singular on 
$\Gamma(t)= \{ \xi(s) ; -\infty <s \leq ct \}$ for each $t\in \mathbb{R}$. 
\end{theorem}

Let us now mention some known results on positive singular solutions of the equation
\begin{equation}\label{eq:porous1}
	u_t=\Delta u^m, \qquad x \in \mathbb{R}^n \setminus \{ \theta(t) \}, \quad  t>0, 
\end{equation}
where $0<m<1$, $n\geq 2$ and $\theta \in C^1([0,\infty); \mathbb{R}^n)$ is a given function. 
We consider \eqref{eq:porous1} with
the initial condition
\begin{equation}\label{eq:initial}
u(x,0)=u_0(x), \qquad x \in \mathbb{R}^n \setminus \{ \theta(0) \}. 
\end{equation}
We are interested in solutions that are singular at $\theta(t)$, that is,
\begin{equation}\label{eq:sing}
  u(x,t) \to \infty \qquad \mbox{ as } x \to \theta(t), \quad t\geq 0.
\end{equation}
For example, when $\theta \equiv 0$ and $n\geq3$, 
\eqref{eq:porous1}  has a singular steady state given by
\[
	\tilde u(x) = C |x|^{-\frac{n-2}{m}},\qquad x\in\mathbb{R}^n\setminus\{0\},
\]
where $C$ is an arbitrary positive constant.  
Another explicit singular solution for $\theta\equiv 0$, $u_0\equiv 0$ and 
\[
m_c:=\frac{(n-2)_+}{n}<m<1
\] 
is
\[
u^*(x,t):=\left(\frac{ct}{|x|^2}\right)^\frac{1}{1-m},\qquad
c:=2m\left(\frac{2}{1-m}-n\right).
\]

For $0<m<1$ and $n\geq 2$, 
one can find in \cite{Ch} a complete classification of nonnegative solutions
of $u_t=\Delta u^m$ in ${\mathcal D}'((\mathbb{R}^n\setminus \{0\})\times(0,\infty))$
which are continuous in $\mathbb{R}^n\times[0,\infty)$ 
with values in $(0,\infty]$, satisfy \eqref{eq:initial} with $u_0\equiv 0$ and \eqref{eq:sing}.

For the existence of self-similar solutions with a standing singularity ($\theta\equiv 0$),
 we refer to \cite{Ch}
when $m_c<m<1$ and to \cite{HuiKim,Vazquez} when $0<m<m_c$ and $n\geq 3$.

If $m_c<m<1$ then all weak solutions of $u_t=\Delta u^m$ with locally
integrable initial data $u_0$ become immediately bounded and continuous,
see \cite{HP}. On the other hand, in the same range $m_c<m<1$, stronger singularities may persist,
see \cite{Chasseigne1, Chasseigne2}. It was shown in \cite{Chasseigne1} that for  $m_c<m<1$,
the strongly
singular set of $u_0$ cannot shrink in time for extended continuous solutions. Here the strongly singular set of $u_0$ is defined as the set of
points at which $u_0$ is not locally integrable 
and an extended continuous solution satisfies the equation pointwise in the
set $\{(x,t)\in\mathbb{R}^{n+1} ;  v(x,t)<\infty\}$, and is continuous with values in $(0,\infty]$. 
The existence of extended continuous solutions with
expanding strongly singular sets is also established in \cite{Chasseigne1}.

The evolution of standing singularities of proper (minimal) solutions of $u_t=\Delta
u^m$ on a bounded domain was studied in \cite{VW} for $n\geq 3$ and
$0<m\leq m_c$.
By a proper solution we mean a solution
obtained as a limit of increasing bounded approximations. 

It was shown in \cite{FTY18} that for $m>(n-2)/(n-1)$ and $n\geq 3$, under some assumptions on given functions $\theta, k, u_0$,
there are solutions of \eqref{eq:porous1}-\eqref{eq:sing} which behave as $k(t)|x-\theta(t)|^{-(n-2)/m}$ 
near $x=\theta(t)$. This corresponds to the singularity of the steady state $\tilde u$. For $m<(n-2)/(n-1)$
such solutions with moving singularities do not exist. The case $n=2$ has been treated recently in \cite{FPTYpre}. 
It has been established there that solutions with moving singularities, 
which behave near the singularity like the fundamental solution of the Laplace equation,
raised to the power $1/m$, exist for $m>0$.

For various results on solutions with moving singularities ($\theta\not\equiv 0$) for the heat equation ($m=1$) we refer to
\cite{KT14,KT16,TY15}, for semilinear heat equations see
\cite{HTY,HY,KT16,KT17,SY09,SY10,SY11,Ta20pre,TY17,TY16} 
and also \cite{KSS20,KZ15} for the Navier-Stokes system.

Our result is somehow disconnected from the previous ones since the nature of the singularity described here is novel and different. Our aim is to contribute to a deeper understanding of singularity formation and non-uniqueness phenomena for the fast diffusion equation. The implications for the Cauchy problem consist in showing the existence of initial functions from which a strong singularity can move 
in time along a prescribed curve, leaving the solution singular behind. As far as we know, this kind of behavior has not been observed 
previously.

The idea of the proof is to use matched asymptotics in order to construct suitable sub- and super-solutions. The most important part of them is derived from an explicit entire solution which can be found in the special case when the curve is a straight line. We also rely on some delicate properties of the distance function. Once entire comparison functions are constructed, the proof can be finished by standard methods.

The paper is organised as follows.
In Section~2 we study an explicit solution when $\{\Gamma(t); t\in \mathbb{R}\}$ is a straight line. In Section~3
we prepare suitable comparison functions for the proof of Theorem~\ref{th:main}. Section~4 contains the
proof and Section~5 a discussion. In 
Appendix~\ref{sec:app}
we give a short derivation of the formula for the explicit
solution from Section~2.

\section{Explicit solution}
In this section, we consider the case $\xi(s)=s\omega$, 
where $\omega \in \mathbb{R}^n$ is a unit vector.

\subsection{Singular traveling wave solution}
Let $c> 0$ be a constant and let $a:=c \omega$ be a velocity vector. 
Set $x=y + ta$. 
Then by taking $v(y,t)=u(y+ta,t)$, we see that $v$ satisfies the equation
\begin{equation}\label{eq:v}
	v_t=\Delta_y v^m + a \cdot \nabla_y v.
\end{equation}
If $m_*< m < 1$, this equation has a stationary solution explicitly expressed as
\begin{align}
	\label{eq:expli}
	&\varphi=\varphi(y) := A ( |a|  |y|+a \cdot y )^{-\frac{1}{1-m}}, \\
	\label{eq:Adef}
	&A=A(n,m):= 
	\left( 
	\frac{(n-1)m}{1-m} \left( m - \frac{n-3}{n-1} \right) 
	\right)^\frac{1}{1-m}. 
\end{align}
We observe that the solution 
\[
	u(x,t)=\varphi(x-ta) 
	= A \left( |a|  |x-ta|+a \cdot (x-ta) \right)^{-\frac{1}{1-m}}
\] 
of \eqref{eq:porous} has a singularity on the set 
$\{s \omega; -\infty< s \leq ct \}$ for each $t\in\mathbb{R}$, 
and so we call this solution the singular traveling wave solution. 

We note that, in the case where $n=1$ and $0<m<1$, 
the explicit solution in \eqref{eq:expli} is known 
as a semi-wavefront solution, cf. for instance \cite[Section 2.1]{GKbook}. 
See Appendix~\ref{sec:app}
for a straightforward derivation of \eqref{eq:expli}-\eqref{eq:Adef}.

\subsection{Stability of the singular traveling wave solution}
By a direct computation, we have 
\[
	\Delta \varphi^m =- a \cdot \nabla  \varphi =
	\frac{A}{1-m} \frac{|a|}{|y|} ( |a||y|+a\cdot y )^{-\frac{1}{1-m}} \geq0. 
\]
Let $0<\gamma<1$. Setting $\varphi^+(y):=(1+\gamma) \varphi(y)$ gives 
\[
\begin{aligned}
	\Delta    (\varphi^+)^m + a \cdot \nabla \varphi^+
	&= (1+\gamma)^m \Delta \varphi^m + (1+\gamma) a \cdot \nabla \varphi \\
	& = (  (1+\gamma)^{m} - (1+\gamma) ) \Delta \varphi^m \leq 0.
\end{aligned}
\]
Hence $\varphi^+$ is a super-solution of \eqref{eq:v}.  Similarly, 
$\varphi^-(y):=(1-\gamma) \varphi(y)$ is a sub-solution of \eqref{eq:v}, 
and so the functions 
\[
	u^\pm (x ,t):=   \varphi^\pm(x-ta)  = (1\pm \gamma) \varphi(x - ta) 
\]
are a super-solution and a sub=solution of \eqref{eq:porous}, respectively. 
Thus, the singular traveling wave solution $\varphi(x-ta)$ is stable, 
and so it is expected that we can construct suitable comparison functions 
for proving Theorem \ref{th:main}.

\subsection{Singular traveling wave solution in cylindrical coordinates}
For $x\in \mathbb{R}^n$, 
let $r(x)$ be the distance between $x$ and the line $\{ \xi(s)=s\omega ; s \in \mathbb{R} \}$. 
Writing
\[
	x=z(x) + s(x)\omega, \qquad z(x) \perp \omega, 
\]
we have $r(x)= |z(x)|=|x-s(x)\omega|$ and $s(x)=\omega \cdot x$. 
Then by using
\[
	x - ta = x - s(x) \omega + s(x) \omega -ct\omega 
	= z(x) +  (s(x)-ct)\omega,
\]
we obtain
\[
\begin{aligned}
	&|a| |x - ta|=c \sqrt{|z(x)|^2 + |(s(x)-ct)\omega|^2} 
	= c \sqrt{r^2(x)+(s(x)-ct)^2}, \\
	&a\cdot (x - ta) = 
	c\omega\cdot z(x) + 
	c\omega \cdot (s(x)-ct)\omega =   c(s(x)-ct).
\end{aligned}
\]
Hence the traveling singular solution can be written as
\[
\begin{aligned}
	&u(x,t)=\psi(r(x),s(x)-ct)
	=A c^{-\frac{1}{1-m}}
	\left( \sqrt{r^2(x)+(s(x)-ct)^2} + s(x)-ct \right)^{-\frac{1}{1-m}}, \\
	&\psi=\psi(r,\sigma)
	:=A c^{-\frac{1}{1-m}} 
	\left( \sqrt{r^2+ \sigma^2} + \sigma \right)^{-\frac{1}{1-m}}.
\end{aligned}
\]
Based on the above observation, 
we handle more general cases in the subsequent sections.

\section{Comparison functions}
In what follows, we take $r_0>0$ such that $0<r_0<\tilde r_0$, 
where $\tilde r_0$ 
 ($<1$) is given in Condition \ref{con:curve}. 
Then for any $x\in \Gamma_{r_0}$, 
there exists a unique real number $s(x)\in\mathbb{R}$ such that 
$r(x)=|x-\xi(s(x))|$. 
We remark that $r(x)$ and $s(x)$ are $C^3$-functions on $\Gamma_{r_0}\setminus\Gamma$, 
since $\xi$ is $C^3$.

The goal of this section is to prove the following proposition,
which guarantees the existence of suitable 
comparison functions for showing Theorem~\ref{th:main}.

\begin{proposition}\label{pro:compari}
Let $n\geq2$, $m_*<m<1$, $c>0$ and $0<\varepsilon<1$. 
Then there exist a super-solution $\overline{u}$ and a sub-solution $\underline{u}$ of \eqref{eq:porous} in $Q$ 
such that the following {\rm(i)} and {\rm(ii)} hold. 
\begin{itemize}
\item[{\rm(i)}]
$\overline{u}\geq \underline{u} >0$ on $Q$. 
\item[{\rm(ii)}]
There exists a constant $0<\delta<r_0$ depending on $c$ and $\varepsilon$ such that 
\[
	(1-\varepsilon) U(x,t)
	\leq \underline{u}(x,t) \leq \overline{u}(x,t)
	\leq (1+\varepsilon) U(x,t) 
\]
for any $(x,t)\in \mathbb{R}^{n+1}$ 
satisfying $0<r(x)\leq \delta$ and $-\infty < s(x)\leq ct+\delta$, 
where $U$ is given in Theorem \ref{th:main}. 
\end{itemize}
\end{proposition}

\subsection{Ingredients of comparison functions}
Our comparison functions are based on the following function. 
\[
\begin{aligned}
	&U(x,t)  = 
	\psi(r(x),\sigma(x,t)) 
	= A c^{-\frac{1}{1-m}}\left( \sqrt{r^2(x)+\sigma^2(x,t)} 
	+ \sigma(x,t) \right)^{-\frac{1}{1-m}}, \\
	&\sigma(x,t) := s(x)-ct, 
\end{aligned}
\]
where $U$ is the same as in Theorem \ref{th:main}.  
Notice that $U$ is defined at least on the set 
\[
	\left\{ (x,t)\in\mathbb{R}^{n+1}; 
	\begin{aligned}
	&\mbox{ either }x\in \Gamma_{\tilde r_0}\setminus\Gamma \mbox{ with }s(x)\leq ct, \\
	&\mbox{ or } x\in \Gamma_{\tilde r_0} \mbox{ with }s(x)> ct 
	\end{aligned}
	\right\}, 
\]
where $\tilde r_0$ is given in Condition \ref{con:curve}. 
We observe that $U(\cdot,t)$ is singular on the set 
$\Gamma(t)=\{ \xi(s); -\infty < s \leq ct \}$ for each $t\in\mathbb{R}$. 
In order to compute the derivatives of $U$, 
we explicitly compute $\nabla s$, $\nabla r$, $\Delta s$ and $\Delta r$ as follows. 
First, we prepare a fundamental lemma.

\begin{lemma}
For $x\in \Gamma_{r_0}\setminus\Gamma$, 
the following equality holds. 
\begin{equation}\label{eq:orth}
	(x-\xi(s(x)))\cdot \xi'(s(x))=0. 
\end{equation}
Moreover, $\xi'(s(x))\cdot \xi''(s(x))=0$ also holds 
for $x\in \Gamma_{r_0}\setminus\Gamma$. 
\end{lemma}

\begin{proof}
Since $\xi(s(x))$ is the nearest point from $x$, 
we have 
\[
\begin{aligned}
	&\partial_s(|x-\xi(s)|^2)|_{s=s(x)} = -2(x-\xi(s(x)))\cdot\xi'(s(x)) =0, \\
	&\partial_s^2(|x-\xi(s)|^2)|_{s=s(x)} 
	= 2 |\xi'(s(x))|^2 - 2(x-\xi(s(x)))\cdot\xi''(s(x)) \geq0. 
\end{aligned}
\]
Then \eqref{eq:orth} follows. Moreover, by $|\xi'(s(x))|^2=1$, 
we also have 
$(x-\xi(s(x)))\cdot \xi''(s(x))\leq 1$. 
From $|\xi'(s)|^2=1$ for $s\in\mathbb{R}$, it follows that 
$\xi'(s)\cdot\xi''(s)=0$ for $s\in\mathbb{R}$. 
\end{proof}

Note that $|(x-\xi(s(x)))\cdot \xi''(s(x))| \leq r_0 K \leq 1/2$ 
by $r_0<\tilde r_0$ and Condition \ref{con:curve}. 
By using this lemma, we can compute $\nabla s$, $\nabla r$, $\Delta s$ and 
$\Delta r$ as follows.

\begin{lemma}\label{lem:deriv}
For $x\in \Gamma_{r_0}\setminus\Gamma$, the following equalities hold. 
\[
\begin{aligned}
	&\nabla s(x) = 
	\frac{\xi'(s(x))}{1- (x-\xi(s(x))) \cdot \xi''(s(x))}, \\
	&\nabla r(x) = \frac{x-\xi(s(x))}{|x-\xi(s(x))|}, \\
	&\Delta s(x) = 
	\frac{(x-\xi(s(x)))\cdot \xi'''(s(x))}{[1-(x-\xi(s(x)))\cdot \xi''(s(x))]^3}, \\
	&\Delta r(x)  = 
	\frac{n-2-(n-1)(x-\xi(s(x))) \cdot \xi''(s(x))}{1- (x-\xi(s(x))) \cdot \xi''(s(x))} 
	r(x)^{-1}. 
\end{aligned}
\]
In particular, $\nabla s(x) \cdot \nabla r(x) =0$. 
\end{lemma}

\begin{proof}
In this proof, we write $s=s(x)$ and $r=r(x)$ for short. 
From the differentiation of \eqref{eq:orth} with respect to $x$, it follows that  
\begin{equation}\label{eq:deriorth}
\begin{aligned}
	0&= 
	\xi'(s) - |\xi'(s)|^2 \nabla s 
	+ \left[ (x-\xi(s)) \cdot \xi''(s) \right] \nabla s \\
	&= 
	\xi'(s) - \nabla s 
	+ \left[ (x-\xi(s)) \cdot \xi''(s) \right] \nabla s.  
\end{aligned}
\end{equation}
Then, 
\[
	\nabla s = 
	\frac{\xi'(s)}{1- (x-\xi(s)) \cdot \xi''(s)}. 
\]
By taking the divergence in \eqref{eq:deriorth}, we have 
\[
\begin{aligned}
	&\xi''(s)\cdot \nabla s - \Delta s 
	+ \xi''(s)\cdot \nabla s 
	- \xi'(s) \cdot \xi''(s) |\nabla s|^2 \\
	&+ (x-\xi(s))\cdot \xi'''(s) |\nabla s|^2 
	+ \left[ (x-\xi(s)) \cdot \xi''(s) \right] \Delta s 
	=0. 
\end{aligned}
\]
This together with $\xi'\cdot \xi''\equiv 0$ shows that 
\[
	- \Delta s + (x-\xi(s))\cdot \xi'''(s) |\nabla s|^2 
	+ \left[ (x-\xi(s)) \cdot \xi''(s) \right] \Delta s 
	=0, 
\]
and so 
\[
	\Delta s = 
	\frac{(x-\xi(s))\cdot \xi'''(s)}{1-(x-\xi(s))\cdot \xi''(s)} 
	|\nabla s|^2 
	= 
	\frac{(x-\xi(s))\cdot \xi'''(s)}{[1-(x-\xi(s))\cdot \xi''(s)]^3}.  
\]

By the differentiation of $r=|x-\xi(s)|$ and by \eqref{eq:orth}, 
we have 
\[
	\nabla r = 
	\frac{x-\xi(s)- [(x-\xi(s))\cdot \xi'(s)]\nabla s }{|x-\xi(s)|} 
	= \frac{x-\xi(s)}{|x-\xi(s)|}. 
\]
By \eqref{eq:orth} again, we also have 
$\nabla s \cdot \nabla r=0$. 
Taking the divergence of $\nabla r = (x-\xi(s))r^{-1}$ yields 
\[
\begin{aligned}
	\Delta r &= 
	n r^{-1} - r^{-1} \xi'(s) \cdot \nabla s 
	- r^{-2} (x-\xi(s))\cdot \nabla r \\
	&= 
	\left[ n-1 - \xi'(s) \cdot \nabla s \right] r^{-1}  \\
	&=  
	\left[ n-1 - \frac{|\xi'(s)|^2}{1- (x-\xi(s)) \cdot \xi''(s)} \right] 
	r^{-1}.  
\end{aligned}
\]
This together with $|\xi'|\equiv1$ shows the desired equality for $\Delta r$. 
\end{proof}

From Lemma \ref{lem:deriv} and Condition \ref{con:curve}, 
the following lemma  immediately follows.

\begin{lemma}\label{lem:derivesti}
For $x\in \Gamma_{r_0}\setminus\Gamma$, the following inequalities hold. 
\[
\begin{aligned}
	&2^{-1} \leq (1+ r_0 K)^{-1} \leq |\nabla s(x)| \leq (1- r_0 K)^{-1} \leq 2, \\
	&-8r_0 K\leq - \frac{r_0 K}{(1- r_0 K)^3}\leq  \Delta s(x)
	\leq \frac{r_0 K}{(1- r_0 K)^3} \leq 8r_0 K,
	\\
	&\frac{n-2 - (n-1)r_0 K}{1+ r_0 K} r(x)^{-1} 
	\leq \Delta r(x) 
	\leq \frac{n-2 + (n-1)r_0 K}{1- r_0 K} r(x)^{-1}.
\end{aligned}
\]
\end{lemma}

We next compute the derivatives of $U$. Set 
\[
	Q_{r_0}:= \left\{ (x,t)\in\mathbb{R}^{n+1}; 
	\begin{aligned}
	&\mbox{ either }x\in \Gamma_{r_0}\setminus\Gamma \mbox{ with }s(x)\leq ct, \\
	&\mbox{ or } x\in \Gamma_{ r_0} \mbox{ with }s(x)> ct 
	\end{aligned}
	\right\}. 
\]

\begin{lemma}\label{lem:derivU}
For $(x,t)\in Q_{r_0}$, 
the following equalities 
hold. 
\[
\begin{aligned}
	&U_t = \frac{1}{1-m} (r^2+\sigma^2)^{-\frac{1}{2}} cU, \\
	&\begin{aligned}
	\Delta U^m 
	&= 
	\left[ \frac{A^{m-1} m}{(1-m)^2} (r^2+\sigma^2)^{-\frac{1}{2}} cU \right]
	\bigg[ -(1-|\nabla s(x)|^2) 
	\left( \frac{(1-m)\sigma^2}{r^2+\sigma^2} 
	+ \frac{\sigma}{\sqrt{r^2+\sigma^2}} \right) \\
	&\quad  
	+(1+m|\nabla s(x)|^2) 
	-(1-m)r \Delta r(x) - (1-m) \left( \sqrt{r^2+\sigma^2} + \sigma \right)\Delta s(x) 
	\bigg], 
	\end{aligned}
\end{aligned}
\]
where $r=r(x)$ and $\sigma=\sigma(x,t)=s(x)-ct$. 
\end{lemma}

\begin{proof}
The time derivative $U_t$ is easy. 
The spatial derivative $\Delta U^m$ is computed as 
\[
	\Delta U^m 
	=  (\psi^m)_{rr}|\nabla r(x)|^2 + (\psi^m)_{\sigma\sigma} |\nabla s(x)|^2 
	+(\psi^m)_r \Delta r(x) + (\psi^m)_\sigma\Delta s(x). 
\]
Straightforward computations show that 
\[
\begin{aligned}
	(\psi^m)_r &= -\frac{A^{m-1} m }{1-m} r(r^2+\sigma^2)^{-\frac{1}{2}} c U, \\
	(\psi^m)_{rr} &= 
	\left[ \frac{A^{m-1} m }{(1-m)^2} 
	\left( 1 - \frac{\sigma}{\sqrt{r^2+\sigma^2}} \right)  
	- \frac{A^{m-1} m }{1-m} \frac{\sigma^2}{ r^2+\sigma^2}  \right]  
	(r^2+\sigma^2)^{-\frac{1}{2}} c U , \\
	(\psi^m)_\sigma &= 
	- \frac{A^{m-1} m }{1-m} \left( \sqrt{r^2+\sigma^2}+\sigma \right) 
	(r^2+\sigma^2)^{-\frac{1}{2}} c U, \\
	(\psi^m)_{\sigma\sigma} &= 
	\left[ 
	\frac{A^{m-1} m }{(1-m)^2} \left( 1+\frac{\sigma}{\sqrt{r^2+\sigma^2}} \right)
	- \frac{A^{m-1} m }{1-m}  \frac{r^2}{r^2+\sigma^2}  
	\right] (r^2+\sigma^2)^{-\frac{1}{2}} c U.   \\
\end{aligned}
\]
From these equalities and $|\nabla r|=1$, the desired equality 
for $\Delta U^m$ follows. 
\end{proof}

\subsection{Super-solution near $\Gamma$}
Fix $0<\varepsilon'<1$. Set 
\[
	u^+(x,t):= 
	(1+\varepsilon') ( U^m (x,t) + 1 )^\frac{1}{m}. 
\]
We check that $u^+$ is a super-solution of \eqref{eq:porous} 
on $Q_{r_0}$ provided that $r_0$ is small enough. 
By Lemma \ref{lem:derivU}, we have 
\begin{equation}\label{eq:vpindert}
\begin{aligned}
	u^+_t &= 
	(1+\varepsilon') \left( 1 + U^{-m} \right)^{\frac{1}{m}-1} U_t\\
	&= 
	\frac{1+\varepsilon'}{1-m} \left( 1 + U^{-m} \right)^{\frac{1-m}{m}} 
	(r^2+\sigma^2)^{-\frac{1}{2}} cU \\
	&= 
	\frac{1+\varepsilon'}{1-m} 
	\left( 1+ A^{-m}  c^{\frac{m}{1-m}}\left( \sqrt{r^2+\sigma^2} 
	+ \sigma \right)^\frac{m}{1-m} 
	\right)^{\frac{1-m}{m}} 
	(r^2+\sigma^2)^{-\frac{1}{2}} cU. 
\end{aligned}
\end{equation}

Since $\Delta (u^+)^m=(1+\varepsilon')^m\Delta U^m$, we will estimate $\Delta U^m$. 
By Lemmas \ref{lem:derivU} and \ref{lem:derivesti}, we have 
\[
\begin{aligned}
	\Delta U^m & \leq 
	\left[ \frac{A^{m-1} m}{(1-m)^2} (r^2+\sigma^2)^{-\frac{1}{2}} cU \right]
	\bigg[ |  1-|\nabla s(x)|^2  | 
	\left( \frac{(1-m)\sigma^2}{r^2+\sigma^2} 
	+ \frac{|\sigma|}{\sqrt{r^2+\sigma^2}} \right)   \\
	&\quad 
	+ 1+ \frac{m}{(1-r_0 K)^2} 
	-(1-m) \frac{n-2 - (n-1)r_0 K}{1+ r_0 K}  \\
	&\quad 
	+8r_0 K (1-m) \left( \sqrt{r^2+\sigma^2} + \sigma \right) 
	\bigg]. 
\end{aligned}
\]
From the estimate of $|\nabla s|$ in Lemma \ref{lem:derivesti}, \eqref{eq:Adef} and 
\begin{equation}\label{eq:sigr2}
	\frac{(1-m)\sigma^2}{r^2+\sigma^2} 
	+ \frac{|\sigma|}{\sqrt{r^2+\sigma^2}} \leq 2, 
\end{equation}
it follows that 
\[
\begin{aligned}
	\Delta U^m
	& \leq 
	\left[ \frac{A^{m-1} m (n-1) }{(1-m)^2} \left( m- \frac{n-3}{n-1} \right)
	\right.\\
	&\quad \left. 
	+\frac{A^{m-1} m}{(1-m)^2} 
	8r_0 K (1-m) \left( \sqrt{r^2+\sigma^2} + \sigma \right) 
	+O(r_0)
	\right]  (r^2+\sigma^2)^{-\frac{1}{2}} cU \\
	&\leq 
	\left[ \frac{ 1 }{1-m} 
	+C r_0 \left( \sqrt{r^2+\sigma^2} + \sigma \right) 
	+O(r_0)
	\right]  (r^2+\sigma^2)^{-\frac{1}{2}} cU
\end{aligned}
\]
with a constant $C=C(m,A,K)>0$ as $r_0\to0$. Thus, 
\[
\begin{aligned}
	\Delta (u^+)^m
	\leq 
	\left[ \frac{ (1+\varepsilon')^m }{1-m} 
	+C r_0  \left( \sqrt{r^2+\sigma^2} + \sigma \right) 
	+O(r_0)
	\right]  (r^2+\sigma^2)^{-\frac{1}{2}} cU
\end{aligned}
\]
with a constant $C=C(m,A,K,\varepsilon')>0$ as $r_0\to0$.

The above computations show that 
\[
\begin{aligned}
	u^+_t  - \Delta (u^+)^m 
	&\geq 
	\frac{(1+\varepsilon')^m}{1-m} 
	\Bigg[ 
	(1+\varepsilon')^{1-m}
	\left( 1+ A^{-m}  c^{\frac{m}{1-m}}\left( \sqrt{r^2+\sigma^2} 
	+ \sigma \right)^\frac{m}{1-m} 
	\right)^{\frac{1-m}{m}}  \\
	&
	\quad 
	- 1 - C r_0  \left( \sqrt{r^2+\sigma^2} + \sigma \right) 
	+O(r_0)
	\Bigg]
	(r^2+\sigma^2)^{-\frac{1}{2}} cU. 
\end{aligned}
\]
We first assume that 
\[
	A^{-m}  c^{\frac{m}{1-m}}\left( \sqrt{r^2+\sigma^2} 
	+ \sigma \right)^\frac{m}{1-m}  \geq2. 
\]
Then, 
\[
\begin{aligned}
	&(1+\varepsilon')^{1-m}
	\left( 1+ A^{-m}  c^{\frac{m}{1-m}}\left( \sqrt{r^2+\sigma^2} 
	+ \sigma \right)^\frac{m}{1-m} 
	\right)^{\frac{1-m}{m}} 
	- 1 - C r_0  \left( \sqrt{r^2+\sigma^2} + \sigma \right) \\
	&\geq 
	(1+\varepsilon')^{1-m}
	 A^{-(1-m)}  c \left( \sqrt{r^2+\sigma^2} + \sigma \right) 
	- 1 - C r_0  \left( \sqrt{r^2+\sigma^2} + \sigma \right) \\
	&\geq 
	2^\frac{1}{m} 
	(1+\varepsilon')^{1-m}- 1
	 + \left( \frac{(1+\varepsilon')^{1-m}}{2} A^{-(1-m)}  c -C r_0 \right)
	\left( \sqrt{r^2+\sigma^2} + \sigma \right) \geq0 
\end{aligned}
\]
provided that $r_0$ is small depending on $\varepsilon'$, $m$, $c$ and $A$. 
On the other hand, we next assume that 
\[
	A^{-m}  c^{\frac{m}{1-m}}\left( \sqrt{r^2+\sigma^2} 
	+ \sigma \right)^\frac{m}{1-m}  \leq2. 
\]
In this case, we have 
\[
\begin{aligned}
	&(1+\varepsilon')^{1-m}
	\left( 1+ A^{-m}  c^{\frac{m}{1-m}}\left( \sqrt{r^2+\sigma^2} 
	+ \sigma \right)^\frac{m}{1-m} 
	\right)^{\frac{1-m}{m}} 
	- 1 - C r_0  \left( \sqrt{r^2+\sigma^2} + \sigma \right) \\
	&\geq 
	(1+\varepsilon')^{1-m} -1 - 2^\frac{1-m}{m} C A^{1-m} c^{-1} r_0
	\geq 0
\end{aligned}
\]
provided that $r_0$ is small depending on $m$, $A$, $c$ and $C=C(m,A,K,\varepsilon')$. 
Hence $u^+$ is a super-solution in $Q_{r_0}$ 
provided that $r_0$ is sufficiently small 
depending only on $m$, $A$, $K$, $c$ and $\varepsilon'$.

\subsection{Super-solution on $Q$}
We construct a super-solution on $Q$, where $Q$ is defined by \eqref{eq:Qdef}. 
For $B, b>1$, define 
\begin{equation}\label{eq:Vpdef}
	\overline{u}(x,t):= 
	\left\{ 
	\begin{aligned}
	& \left[ \eta(r(x)) 
	(u^+)^m (x,t) + B(b-\eta(r(x))) \right]^\frac{1}{m} 
	&&\mbox{ if }r(x)\leq r_0, \\
	& (B b)^\frac{1}{m} 
	&&\mbox{ otherwise}. 
	\end{aligned}
	\right.
\end{equation}
Here $\eta$ is defined below. 
Let $r_1$ and $r_2$ satisfy $0<r_1<r_2<r_0$. 
We define $\eta\in C^\infty([0,\infty))$ by 
\[
	\eta=\eta(\rho) := 
	\left\{ 
	\begin{aligned}
		& 1 && \mbox { for } 0\leq \rho \leq r_1, \\
		& 
		\frac{ \Phi(\rho) }{ \Phi(\rho) + \Psi(\rho)}
		&&\mbox{ for }r_1< \rho < r_2, \\
		& 0 && \mbox { for } \rho \geq r_2, 
	\end{aligned}
	\right.
\]
where 
\[
	\Phi(\rho):= \exp\left( -\frac{1}{r_2-\rho} \right), \qquad 
	\Psi(\rho):= \exp\left( -\frac{1}{\rho-r_1} \right). 
\]
After straightforward computations, we have 
\[
\begin{aligned}
	&\Phi'=-\frac{1}{(r_2-\rho)^2}\Phi, \qquad 
	\Phi''=\frac{1-2r_2+2\rho}{(r_2-\rho)^4}\Phi, \\
	&\Psi'=\frac{1}{(\rho-r_1)^2}\Psi, \qquad 
	\Psi''=\frac{1+2r_1-2\rho}{(\rho-r_1)^4}\Psi, \qquad 
	\eta' = \frac{\Phi'\Psi-\Phi\Psi'}{(\Phi+\Psi)^2}, \\
	&\eta''=\frac{1}{(\Phi+\Psi)^2} (\Phi''\Psi-\Phi\Psi'')
	-\frac{2}{(\Phi+\Psi)^3} (\Phi'\Psi-\Phi\Psi')(\Phi'+\Psi'), 
\end{aligned}
\]
and so 
\[
\begin{aligned}
	\eta'(\rho) &= -\frac{\Phi\Psi}{(\Phi+\Psi)^2} 
	\left( \frac{1}{(r_2-\rho)^2} + \frac{1}{(\rho-r_1)^2} \right), \\
	\eta''(\rho) 
	&= \frac{\Phi\Psi}{(\Phi+\Psi)^2} 
	\left( 
	\frac{1-2r_2+2\rho}{(r_2-\rho)^4} 
	-\frac{1+2r_1-2\rho}{(\rho-r_1)^4} 
	\right) \\
	&\quad 
	- \frac{2 \Phi\Psi}{(\Phi+\Psi)^3}
	\left( \frac{1}{(r_2-\rho)^2} + \frac{1}{(\rho-r_1)^2} \right) 
	\left( \frac{\Phi}{(r_2-\rho)^2} - \frac{\Psi}{(\rho-r_1)^2} \right)
\end{aligned}
\]
for $r_1< \rho< r_2$. 
We note that
\begin{equation}\label{eq:limr2}
	\lim_{\rho\uparrow r_2} \frac{\eta(\rho)}{\eta''(\rho)} 
	= \lim_{\rho\uparrow r_2} \frac{|\eta'(\rho)|}{\eta''(\rho)} 
	=0. 
\end{equation}

By direct computations, we have 
\[
\begin{aligned}
	&\overline{u}_t =
	[(u^+)^m\eta  + B(b-\eta)]^{\frac{1-m}{m}}  (u^+)^{m-1} 
	\eta u^+_t, \\
	&
	\begin{aligned}
	-\Delta \overline{u}^m
	&= (B-(u^+)^m) \Delta(\eta(r(x))) 
	- 2\nabla (\eta(r(x))) \cdot \nabla (u^+)^m 
	- \eta(r(x))\Delta(u^+)^m \\
	&= (B-(u^+)^m) (\eta''+\eta'\Delta r) 
	- 2 \eta' \nabla r \cdot \nabla (u^+)^m 
	- \eta \Delta(u^+)^m, 
	\end{aligned}
\end{aligned}
\]
where $\nabla r$ and $\Delta r$ are evaluated at $r(x)$. 
In the region $r(x)\leq r_1$, we have 
\[
\begin{aligned}
	\overline{u}_t-\Delta \overline{u}^m &=
	[(u^+)^m  + B(b-1)]^{\frac{1-m}{m}}  (u^+)^{m-1} u^+_t 
	- \Delta(u^+)^m 
	\geq 
	u^+_t - \Delta(u^+)^m \geq0, 
\end{aligned}
\]
since $u^+_t>0$ by \eqref{eq:vpindert}. 
On the other hand, in the region $r(x)\geq r_2$, we have 
$\overline{u}_t-\Delta \overline{u}^m = 0$. 
Hence it suffices to consider the region $r_1< r(x) < r_2$.

First, we observe the case where $r(x)$ is smaller than $r_2$ and is close to $r_2$. 
By $u^+_t>0$, we have $\overline{u}_t\geq0$, and so 
\[
	\overline{u}_t -\Delta \overline{u}^m \geq  
	(B-(u^+)^m) (\eta''+\eta'\Delta r)  
	- 2 \eta' \nabla r \cdot \nabla (u^+)^m  
	- \eta \Delta(u^+)^m.  
\]
We take constants  $\tilde C>1$  and $B$ such that 
$\tilde C^{-1} \leq u^+\leq \tilde C$ for $r_1<r(x)<r_2$ and $B>\tilde C^m$. 
By \eqref{eq:limr2} and $\eta''>0$ near $r_2$, we have 
$\eta''-|\eta'||\Delta r|>0$ near $r_2$. 
Then there exists a constant $C>1$ independent of $b$ such that 
\[
	\overline{u}_t -\Delta \overline{u}^m \geq 
	(B-\tilde C^m) (\eta''-C |\eta'|) - C |\eta'|  - C \eta. 
\]
Then by \eqref{eq:limr2} again, there exists a constant $r_1<r_2'<r_2$  
independent of $b$ such that 
$\overline{u}_t -\Delta (\overline{u})^m \geq 0$ for $r_2' <  r(x) < r_2$.

We next examine the case where $r_1< r(x)< r_2'$. 
Note that there exists a constant $\tilde c>0$ such that 
$\tilde c\leq \eta(r(x)) \leq 1$ for $r_1< r(x)< r_2'$. 
By $u^+_t>0$, 
there exists a constant $C>1$ depending on $r_1$ and $r_2$ 
but not on $b$ such that 
\[
\begin{aligned}
	&\overline{u}_t 
	\geq 
	B^{\frac{1-m}{m}} (b-\eta(r(x)) )^{\frac{1-m}{m}}  (u^+)^{m-1} 
	\eta u^+_t 
	\geq 
	\tilde c C^{-1}  B^\frac{1-m}{m} (b-1)^\frac{1-m}{m},  \\
	&-\Delta \overline{u}^m
	=  - \eta \Delta(u^+)^m
	+ (B-(u^+)^m) (\eta''+\eta'\Delta r) 
	- 2 \eta' \nabla r \cdot \nabla (u^+)^m 
	\geq  
	- C - B C. 
\end{aligned}
\]
Then there exists $b>1$ such that 
$\overline{u}_t -\Delta \overline{u}^m >0$ for $r_1 <  r(x) < r_2'$. 
Hence $\overline{u}$ is a super-solution on $Q$.

\subsection{Sub-solution on $Q$}
For $M>0$, set 
\begin{equation}\label{eq:vmdef}
	u^-(x,t):= 
	\left\{ 
	\begin{aligned}
	&(1-\varepsilon') \left[ U^m (x,t) - M 
	- M|\sigma(x,t)|^\frac{m}{1-m} \zeta(\sigma(x,t)) \right]_+^\frac{1}{m}
	&&\mbox{ if }r(x)\leq  r_0, \\
	&0 &&\mbox{ otherwise}, 
	\end{aligned}
	\right. 
\end{equation}
where $\sigma(x,t)=s(x)-ct$, 
$[\cdot]_+$ is the positive part and 
$\zeta\in C^\infty(\mathbb{R})$ is a decreasing function satisfying 
$\zeta(\sigma)=1$ if $\sigma\leq -2$, 
$\zeta(\sigma)=0$ if $\sigma\geq -1$ and $0\leq \zeta \leq 1$. 
Let $M$ satisfy 
\begin{equation}\label{eq:asM}
	M> \max\left\{ 
	3^\frac{m}{1-m}  A^m c^{-\frac{m}{1-m}}  r_0^{-\frac{2m}{1-m}}, 
	10^\frac{m}{1-m}  A^m c^{-\frac{m}{1-m}} r_0^{-\frac{2m}{1-m}}  \right\}. 
\end{equation}

We will see that $u^-$ is a sub-solution of \eqref{eq:porous} on $Q$. 
By the fact that the maximum of two sub-solutions is also a sub-solution and 
Lemmas \ref{lem:negar} and \ref{lem:negs} below, 
we only have to consider the case where 
\begin{equation}\label{eq:Qsas}
	(x,t) \in Q_{r_0}\cap \{u^->0\} 
	\quad \mbox{ and }\quad \sigma(x,t) \leq 1. 
\end{equation}

\begin{lemma}\label{lem:negar}
If $r(x)=r_0$, then $u^-(x,t)=0$.  
\end{lemma}

\begin{proof}
Let $x$ satisfy $r(x)=r_0$. Then, 
\begin{equation}\label{eq:Uesnega}
\begin{aligned}
	U(x,t) & = Ac^{-\frac{1}{1-m}} 
	\left( \sqrt{r_0^2 + \sigma^2} + \sigma \right)^{-\frac{1}{1-m}} \\
	& = 
	Ac^{-\frac{1}{1-m}} |\sigma|^{-\frac{1}{1-m}} 
	\left( \sqrt{1+ (r_0/\sigma)^2} + \sigma/|\sigma| \right)^{-\frac{1}{1-m}} \\
	&\leq 
	Ac^{-\frac{1}{1-m}} |\sigma|^{-\frac{1}{1-m}} 
	\left( \sqrt{1+ (r_0/\sigma)^2} -1 \right)^{-\frac{1}{1-m}}. 
\end{aligned}
\end{equation}

We first consider the case $\sigma\leq -2$. 
From $r_0\leq 1$ and $\sigma^2\geq4$, it follows that 
\[
\begin{aligned}
	\sqrt{1+ ( r_0/\sigma)^2} -1 & = 
	\left( \frac{r_0}{\sigma} \right)^2 
	\int_0^1 \frac{1}{2} 
	\left( 1+ \left( \frac{r_0}{\sigma} \right)^2 \theta \right)^{-\frac{1}{2}} 
	d\theta \\
	&\geq 
	\frac{1}{2}  \left( \frac{r_0}{\sigma} \right)^2 
	\left( 1+ \left( \frac{r_0}{\sigma} \right)^2 \right)^{-\frac{1}{2}} 
	\geq \frac{1}{3}  \left( \frac{r_0}{\sigma} \right)^2, 
\end{aligned}
\]
so that 
\[
\begin{aligned}
	U(x,t) & \leq 
	3^\frac{1}{1-m}  Ac^{-\frac{1}{1-m}}  r_0^{-\frac{2}{1-m}} 
	|\sigma|^{\frac{1}{1-m}}. 
\end{aligned}
\]
This together with \eqref{eq:asM} and $\zeta(\sigma)=1$ for $\sigma\leq -2$ gives 
\[
	\left[ 
	U^m - M 
	- M|\sigma|^\frac{m}{1-m} \zeta(\sigma)  \right]_+
	\leq 
	\left[ 
	\left( 3^\frac{m}{1-m}  A^m c^{-\frac{m}{1-m}}  r_0^{-\frac{2m}{1-m}} 
	- M\right) |\sigma|^\frac{m}{1-m} \right]_+
	= 0. 
\]

We next consider the case $-2\leq \sigma<  0$. 
From $1\leq 16/\sigma^2$ and $r_0<1$,  it follows that 
\[
\begin{aligned}
	\sqrt{1+ ( r_0/\sigma)^2} -1 
	\geq 
	\frac{1}{2}  \left( \frac{r_0}{\sigma} \right)^2 
	\left( 1+ \left( \frac{r_0}{\sigma} \right)^2 \right)^{-\frac{1}{2}} \geq 
	\frac{1}{2}  \left( \frac{r_0}{\sigma} \right)^2 
	\left( \frac{16+ r_0^2}{\sigma^2}  \right)^{-\frac{1}{2}}
	\geq 
	\frac{r_0^2}{10|\sigma|}. 
\end{aligned}
\]
By \eqref{eq:Uesnega}, we have 
\[
\begin{aligned}
	U(x,t) &\leq 
	10^\frac{1}{1-m}  Ac^{-\frac{1}{1-m}} r_0^{-\frac{2}{1-m}}, 
\end{aligned}
\]
and so 
\[
	\left[ 
	U^m - M 
	- M|\sigma|^\frac{m}{1-m} \zeta(\sigma)  \right]_+
	\leq 
	\left[ 10^\frac{m}{1-m}  A^m c^{-\frac{m}{1-m}} r_0^{-\frac{2m}{1-m}}  
	- M \right]_+ = 0. 
\]

Finally, we examine the case $\sigma\geq0$. In this case, 
\[
	\left[ 
	U^m - M 
	- M|\sigma|^\frac{m}{1-m} \zeta(\sigma)  \right]_+
	\leq 
	\left[ A^m c^{-\frac{m}{1-m}} r_0^{-\frac{m}{1-m}} 
	- M \right]_+ = 0. 
\]
The lemma follows. 
\end{proof}

\begin{lemma}\label{lem:negs}
If $\sigma(x,t)\geq1$, then $u^-(x,t)=0$. 
\end{lemma}

\begin{proof}
If $\sigma\geq1$, then 
\[
\begin{aligned}
	&U(x,t) = Ac^{-\frac{1}{1-m}} 
	\left( \sqrt{r^2(x) + \sigma^2} + \sigma \right)^{-\frac{1}{1-m}} 
	\leq  A  c^{-\frac{1}{1-m}}, \\
	&\left[ U^m - M 
	- M|\sigma(x,t)|^\frac{m}{1-m} \zeta(\sigma(x,t))  \right]_+
	\leq 
	\left[ A^m  c^{-\frac{m}{1-m}} - M \right]_+ = 0, 
\end{aligned}
\]
the lemma follows. 
\end{proof}

We consider the case \eqref{eq:Qsas} with 
$\sigma\leq -2$, $-2\leq \sigma\leq -1$ and $-1 \leq \sigma \leq 1$, respectively. 
First, we assume \eqref{eq:Qsas} with $\sigma\leq -2$. 
By the negativity of $\sigma$ and the positivity of $U_t$, we have 
\[
\begin{aligned}
	u^-_t &= 
	(1-\varepsilon') \frac{1}{m} 
	\left[  U^m - M - M |\sigma|^\frac{m}{1-m} \right]^{\frac{1}{m}-1} 
	\left[ 
	m U^{m-1} U_t + \frac{Mm}{1-m} c |\sigma|^{\frac{m}{1-m}-2} \sigma
	\right] \\
	&\leq 
	(1-\varepsilon') \frac{1}{m} U^{1-m} \left[ 
	m U^{m-1} U_t + \frac{Mm}{1-m} c |\sigma|^{\frac{m}{1-m}-2} \sigma
	\right] \\
	&\leq 
	(1-\varepsilon') U_t
	= \frac{1-\varepsilon'}{1-m} (r^2+\sigma^2)^{-\frac{1}{2}} cU. 
\end{aligned}
\]
Direct computations yield 
\[
	\Delta(u^-)^m 
	= (1-\varepsilon')^m \Delta U^m - (1-\varepsilon')^m M\Delta(|\sigma|^\frac{m}{1-m}). 
\]
By Lemmas \ref{lem:derivU} and \ref{lem:derivesti}, we have 
\[
\begin{aligned}
	\Delta U^m 
	&\geq 
	\left[ \frac{A^{m-1} m}{(1-m)^2} (r^2+\sigma^2)^{-\frac{1}{2}} cU \right]
	\bigg[ -| 1-|\nabla s|^2|  
	\left( \frac{(1-m)\sigma^2}{r^2+\sigma^2} 
	+ \frac{|\sigma|}{\sqrt{r^2+\sigma^2}} \right)  \\
	&\quad  
	+ 1+ \frac{m}{(1+r_0K)^2} 
	-(1-m) \frac{n-2+(n-1)r_0K}{1-r_0K}  \\
	&\quad 
	- 8r_0 K (1-m) \left( \sqrt{r^2+\sigma^2} + \sigma \right) 
	\bigg]. 
\end{aligned}
\]
From the negativity of $\sigma$, it follows that 
\begin{equation}\label{eq:rsig4}
\begin{aligned}
	\sqrt{r^2+\sigma^2} + \sigma 
	&= |\sigma| \left( \sqrt{ 1+ (r/\sigma)^2 } -1 \right)  \\
	&=\frac{r^2}{|\sigma|} 
	\int_0^1 \frac{1}{2} 
	\left( 1+ \left( \frac{r}{|\sigma|} \right)^2 \theta \right)^{-\frac{1}{2}} d\theta 
	\leq 
	\frac{r^2}{2|\sigma|}. 
\end{aligned}
\end{equation}
By $|\sigma|\geq 2$ and $r_0\leq 1$, we have 
\[
	\sqrt{r^2+\sigma^2} + \sigma 
	\leq \frac{1}{4}. 
\]
This together with \eqref{eq:sigr2} 
implies that 
\[
\begin{aligned}
	\Delta U^m 
	&\geq 
	\left[ \frac{A^{m-1} m}{(1-m)^2} (r^2+\sigma^2)^{-\frac{1}{2}} cU \right]
	\left[ -2 \left( 1-\frac{1}{(1+r_0K)^2} \right)   \right.\\
	&\quad \left. 
	+1+\frac{m}{(1+r_0K)^2} -(1-m) \frac{n-2+(n-1)r_0K}{1-r_0K}  - 2r_0 K (1-m) 
	\right] \\
	&=
	\left[ \frac{A^{m-1} m(n-1)}{(1-m)^2} \left( m-\frac{n-3}{n-1} \right) + O(r_0) \right]
	(r^2+\sigma^2)^{-\frac{1}{2}} cU 
	\qquad (r_0\to0). 
\end{aligned}
\]

We estimate $\Delta (|\sigma|^\frac{m}{1-m})$. 
By computations, we have 
\[
	\Delta (|\sigma|^\frac{m}{1-m}) 
	= \frac{m}{1-m}\left( \frac{m}{1-m}-1 \right) 
	|\sigma|^{\frac{m}{1-m}-2}|\nabla s|^2 
	+ \frac{m}{1-m} |\sigma|^{\frac{m}{1-m}-2} \sigma \Delta s, 
\]
and so 
\[
\begin{aligned}
	| \Delta (|\sigma|^\frac{m}{1-m}) |
	&\leq 
	C  \left[ 
	(r^2+\sigma^2)^{\frac{1}{2}} c^{-1} U^{-1} 
	\left( 
	|\sigma|^{\frac{m}{1-m}-2} |\nabla s|^2 
	+ |\sigma|^{\frac{m}{1-m}-1} |\Delta s|
	\right) \right] \\
	&\quad 
	\times (r^2+\sigma^2)^{-\frac{1}{2}} cU
\end{aligned}
\]
with a constant $C=C(m)>0$. 
Then, \eqref{eq:rsig4}, $r^2\leq \sigma^2$, Lemma \ref{lem:derivesti} 
and $|\sigma|\geq 2$ show that 
\[
\begin{aligned}
	&(r^2+\sigma^2)^{\frac{1}{2}} c^{-1} U^{-1} 
	\left( |\sigma|^{\frac{m}{1-m}-2} |\nabla s|^2 
	+ |\sigma|^{\frac{m}{1-m}-1} |\Delta s| \right) \\
	&= (r^2+\sigma^2)^{\frac{1}{2}} 
	A^{-1} c^{\frac{1}{1-m}-1} \left( \sqrt{r^2+\sigma^2} + \sigma \right)^\frac{1}{1-m}
	\left( |\sigma|^{\frac{m}{1-m}-2} |\nabla s|^2 
	+ |\sigma|^{\frac{m}{1-m}-1} |\Delta s| \right) \\
	&\leq 
	C|\sigma| \left( \frac{r_0^2}{|\sigma|} \right)^\frac{1}{1-m} 
	\left( 
	4 |\sigma|^{\frac{m}{1-m}-2} 
	+ 8r_0 K |\sigma|^{\frac{m}{1-m}-1} 
	\right) \\
	&\leq 
	C r_0^\frac{2}{1-m} 
	\left(  4 |\sigma|^{-2}  + 8 K |\sigma|^{-1}  \right) 
	\leq 
	C r_0^\frac{2}{1-m} 
	\left(  1  + 4 K \right) 
	=O(r_0^\frac{2}{1-m} ) \qquad (r_0\to0) 
\end{aligned}
\]
with some constant $C=C(m,A,c)>0$. 

From the above computations and \eqref{eq:Adef}, it follows that 
\begin{equation}\label{eq:vinmasym}
\begin{aligned}
	&u^-_t -\Delta (u^-)^m \\
	&\leq 
	\left[ \frac{1-\varepsilon'}{1-m} - (1-\varepsilon')^m 
	\frac{A^{m-1} m(n-1)}{(1-m)^2} \left( m-\frac{n-3}{n-1} \right) + o(1)
	\right] \\
	&\quad \times (r^2+\sigma^2)^{-\frac{1}{2}} cU \\
	&= 
	- (1-\varepsilon')^m \left[ 
	\frac{1- (1-\varepsilon')^{1-m}}{1-m}   + o(1) 
	\right](r^2+\sigma^2)^{-\frac{1}{2}} cU
	\qquad (r_0\to0). 
\end{aligned}
\end{equation}
Hence $u^-$ is a sub-solution in the case \eqref{eq:Qsas} with $\sigma\leq -2$ 
provided that $r_0$ is sufficiently small. 
Remark that the smallness of $r_0$ is determined only by 
$m$, $n$, $A$, $c$, $K$, $M$ and $\varepsilon'$.

Let us next consider the case \eqref{eq:Qsas} with $-2\leq \sigma\leq -1$. 
By the negativity of $\sigma$ and $\zeta'$, we have 
\[
\begin{aligned}
	u^-_t &= 
	(1-\varepsilon') \frac{1}{m} 
	\left[  U^m - M - M |\sigma|^\frac{m}{1-m} \zeta(\sigma) \right]^{\frac{1}{m}-1} \\
	&\quad \times 
	\left[ 
	m U^{m-1} U_t + \frac{Mm}{1-m} c |\sigma|^{\frac{m}{1-m}-2} \sigma \zeta(\sigma)
	+ M|\sigma|^\frac{m}{1-m} c \zeta'(\sigma) 
	\right] 
	\leq 
	(1-\varepsilon') U_t. 
\end{aligned}
\]
From the smoothness of $|\sigma|^\frac{m}{1-m}\zeta(\sigma)$ 
as a function for $-1\leq \sigma \leq -2$, it follows that 
\[
\begin{aligned}
	\Delta(u^-)^m 
	&= (1-\varepsilon')^m \Delta U^m - (1-\varepsilon')^m M\Delta(|\sigma|^\frac{m}{1-m} \zeta(\sigma) ) \\
	&\geq (1-\varepsilon')^m \Delta U^m - C \left[ 
	(r^2+\sigma^2)^{\frac{1}{2}} c^{-1} U^{-1} \right] (r^2+\sigma^2)^{-\frac{1}{2}} cU
\end{aligned}
\]
with a constant $C=C(\varepsilon', M, \zeta)>0$. 
Since \eqref{eq:rsig4} also holds for the case $-2\leq \sigma\leq -1$, 
we see from $r^2 +\sigma^2 \leq 5$ and $|\sigma|\geq 1$ that 
\[
\begin{aligned}
	(r^2+\sigma^2)^{\frac{1}{2}} c^{-1} U^{-1} 
	&=
	(r^2+\sigma^2)^{\frac{1}{2}} c^{\frac{1}{1-m}-1} A^{-1} 
	\left( \sqrt{r^2+\sigma^2} + \sigma \right)^\frac{1}{1-m} \\
	&\leq 
	C \left( \sqrt{r^2+\sigma^2} + \sigma \right)^\frac{1}{1-m} 
	\leq 
	C \left( \frac{r_0^2}{|\sigma|} \right)^\frac{1}{1-m} 
	\leq C r_0^\frac{2}{1-m} = O(r_0^\frac{2}{1-m}) 
\end{aligned}
\]
as $r_0\to0$, where $C=C(m,A,c)>0$ is some constant. 
Hence we obtain the same estimates
as in \eqref{eq:vinmasym}. 
Then we see that 
$u^-$ is a sub-solution in the case \eqref{eq:Qsas} with $-2\leq \sigma\leq -1$ 
if $r_0$ is sufficiently small depending only on 
$m$, $n$, $A$, $c$, $K$, $M$, $\varepsilon'$ and $\zeta$.

Finally, we examine the case \eqref{eq:Qsas} with $-1 \leq \sigma\leq 1$. 
In this case, we have 
\[
\begin{aligned}
	&u^-_t =  (1-\varepsilon') \frac{1}{m} 
	\left[  U^m - M \right]^{\frac{1}{m}-1}  m U^{m-1} U_t 
	\leq 
	(1-\varepsilon') U_t, \\
	& \Delta(u^-)^m = (1-\varepsilon')^m \Delta U^m. 
\end{aligned}
\]
Then the same estimates
as in \eqref{eq:vinmasym} immediately follow, and so 
$u^-$ is a sub-solution in the case \eqref{eq:Qsas} with $-1\leq \sigma\leq 1$ 
if $r_0$ is sufficiently small depending only on 
$m$, $n$, $A$, $K$ and $\varepsilon'$.

Recall that we only have to consider the case \eqref{eq:Qsas}. 
Thus, $u^-$ is a sub-solution of \eqref{eq:porous} on $Q$ 
provided that $r_0$ is sufficiently small depending only on 
$m$, $n$, $A$, $c$, $K$, $M$, $\varepsilon'$ and $\zeta$.

\subsection{Positive comparison functions}
We prove Proposition \ref{pro:compari}.

\begin{proof}[Proof of Proposition Proposition \ref{pro:compari}] 
Let $0<\varepsilon'<\varepsilon<1$. 
We set $\overline{u}$ and $u^-$ as in \eqref{eq:Vpdef} and \eqref{eq:vmdef}, respectively. 
Define
\[
	\underline{u}(x,t):= \max\{ u^-(x,t), \varepsilon \}. 
\]
Since the maximum of two sub-solutions is also a sub-solution, 
$\underline{u}$ is a positive sub-solution on $Q$. 
Moreover, we can easily check that $\underline{u}\leq \overline{u}$ on $Q$. Then (i) holds. 

We prove (ii). 
By the choice of $\overline{u}$ and $\underline{u}$, we have 
\[
\begin{aligned}
	&\overline{u} = 
	\left[ (1+\varepsilon')^m ( U^m + 1 ) + B(b-1) \right]^\frac{1}{m}, \\
	&\underline{u} = 
	(1-\varepsilon') \left[ U^m - M 
	- M|\sigma|^\frac{m}{1-m} \zeta(\sigma) \right]_+^\frac{1}{m} 
\end{aligned}
\]
for $(x,t)\in \Gamma_{r_1}\setminus \Gamma$. 
For $(x,t)\in\mathbb{R}^{n+1}$ with 
$0<r(x)\leq \delta$ and $\sigma \leq \delta$, by \eqref{eq:rsig4}, we have 
\[
\begin{aligned}
	U^m &= 
	A^m c^{-\frac{m}{1-m}} 
	\left( \sqrt{r^2+\sigma^2}+\sigma \right)^{-\frac{m}{1-m}} \\
	&\geq 
	\left\{ 
	\begin{aligned}
	&\left( \sqrt{2}+1 \right)^{-\frac{m}{1-m}}  
	A^m c^{-\frac{m}{1-m}} \delta^{-\frac{m}{1-m}} 
	&&\mbox{ if }-\infty <\sigma \leq \delta,  \\
	&A^m c^{-\frac{m}{1-m}} 
	\left( \frac{\delta^2}{2|\sigma|}  \right)^{-\frac{m}{1-m}} 
	&&\mbox{ if }-\infty< \sigma<0. 
	\end{aligned} 
	\right.
\end{aligned}
\]
Then there exists a constant $C=C(m,A,c)>0$ such that 
\[
	1\leq  
	\left\{ 
	\begin{aligned}
	&C \delta^{\frac{m}{1-m}} U^m
	&&\mbox{ if }-\infty <\sigma \leq \delta,  \\
	&C |\sigma|^{-\frac{m}{1-m}}
	\delta^\frac{2m}{1-m} U^m
	&&\mbox{ if }-\infty<\sigma<0. 
	\end{aligned}
	\right.
\]
From this, it follows that 
\[
	\overline{u}  \leq 
	\left[ (1+\varepsilon')^m ( 1 + C\delta^\frac{m}{1-m} ) 
	+ C B(b-1) \delta^\frac{m}{1-m} \right]^\frac{1}{m} U
	\leq (1+\varepsilon) U
\]
and 
\[
\begin{aligned}
	\underline{u} & =
	\left\{ 
	\begin{aligned}
	& (1-\varepsilon') \left[ U^m - M 
	- M|\sigma|^\frac{m}{1-m} \zeta(\sigma) \right]_+^\frac{1}{m} 
	&&\qquad (-\infty<\sigma<0) \\
	&(1-\varepsilon') \left[ U^m - M \right]_+^\frac{1}{m} 
	&&\qquad (0\leq \sigma \leq \delta) \\
	\end{aligned}
	\right. \\
	&\geq 
	\left\{ 
	\begin{aligned}
	& (1-\varepsilon') \left[ 1 - C M \delta^\frac{m}{1-m} 
	- C M \delta^\frac{2m}{1-m} 
	\right]_+^\frac{1}{m}  U
	&&\qquad (-\infty<\sigma<0) \\
	&(1-\varepsilon') \left[ 1 - M C \delta^{\frac{m}{1-m}} \right]_+^\frac{1}{m} U
	&&\qquad (0\leq \sigma \leq \delta) \\
	\end{aligned}
	\right.  \\
	&\geq (1-\varepsilon) U
\end{aligned}
\] 
for $(x,t)\in\mathbb{R}^{n+1}$ with 
$0<r(x)\leq \delta$ and $-\infty< \sigma \leq \delta$ 
provided that $\delta$ is sufficiently small 
depending only on $m$, $A$, $b$, $B$, $c$, $M$, $\varepsilon$ and $\varepsilon'$. 
Hence (ii) follows. 
\end{proof}

\section{Proof of the main theorem}
We define an exhaustion of $Q$, 
and then we apply an argument from
\cite{FMN04} (see also \cite[Lemma 2.1]{CG05}) 
to show the existence of an entire-in-time singular solution.

\begin{proof}[Proof of Theorem \ref{th:main}]
Let $\{ \Omega_i(t)\}_{t\in\mathbb{R}}$ be a family of smooth bounded domains in $\mathbb{R}^n$ 
such that 
$\Omega_i(t)\subsetneq \Omega_{i+1}(t)$, 
$\bigcup_{i\geq 1}\Omega_i(t) = \mathbb{R}^n\setminus\Gamma(t)$ for each $t\in\mathbb{R}$ 
and $S_i:=\{(x,t)\in\mathbb{R}^{n+1}; x\in \partial \Omega_i(t), t\in (-i,i)\}$ 
is smooth. Define 
\[
	Q_i:=\{(x,t)\in\mathbb{R}^n; x\in \Omega_i(t), t\in(-i,i)\}. 
\]
Note that $\bigcup_{i\geq 1}Q_i = Q$ and $Q_i\subsetneq Q_{i+1}$. 
Consider the following approximate problem. 
\[
\left\{ 
\begin{aligned}
	&w_t =\Delta w^m &&\mbox{ in }Q_i, \\
	&w=\underline{u} &&\mbox{ on }S_i, \\
	&w(\cdot,-i)=\underline{u}(\cdot,-i)  &&\mbox{ in }\Omega_i(-i). 
\end{aligned}
\right. 
\]
By the assumptions on $\Omega_i$ and 
the uniform positivity of comparison functions $\overline{u}$ and 
$\underline{u}$ on $Q_i$, 
this approximate problem has a bounded solution $w_i$ satisfying 
$\underline{u}\leq w_i\leq \overline{u}$ in $Q_i$. 
Since $\underline{u}$ is a sub-solution, 
the comparison principle for bounded solutions implies that 
$w_i(\cdot,-i) = \underline{u}(\cdot,-i)
\leq w_{i+1}(\cdot,-i)$ in $\Omega_i(-i)$. 
From the comparison principle for bounded solutions again, it follows that 
\[
	\underline{u}\leq w_i \leq w_{i+1} \leq \overline{u}
	\qquad \mbox{ in } Q_i 
\]
for each $i$. 
Hence the limiting function 
\[
	u(x,t) := \lim_{i\to\infty} w_i(x,t), \qquad 
	(x,t) \in Q
\]
exists and satisfies $\underline{u}\leq u\leq \overline{u}$ in $Q$. 
By the same argument as in \cite[Lemma 5.1]{FTY18} based 
on the parabolic interior regularity theory and a diagonalization argument, 
we see that $u\in C^{2,1}(Q)$ and $w_i\to u$ in $C^{2,1}_{{\rm loc}}(Q)$ as $i\to\infty$. 
Hence $u$ satisfies \eqref{eq:porous} in $Q$. 
Moreover, the desired estimate on $u$ 
immediately follows from Proposition \ref{pro:compari}, 
and the proof is complete. 
\end{proof}

\section{Discussion}
While our focus here is on furthering the classification of singular behavior in fast nonlinear diffusion, rather than
upon applications, we now briefly comment on the nature of latter that lies in the background. 
A key phenomenon associated with fast diffusion is the suppression of transport at high ``concentrations" $u$ 
(see \cite{King} for a number of illustrative applications).
It is hoped that the associated intuition clarifies the physical status of singular solutions in the context of localized
sources of material, the line singularities in the above being associated with restricted ability of material to diffuse away
from the ridge of high concentration laid down by the moving source at the head of the snake (possibly augmented by
continued injection along the evolving line).

The snaking solutions are the simplest representatives of a much more general class, that is in turn illustrative of very
wide-ranging issues of non-uniqueness in the equation of fast diffusion (cf. \cite{FTY18} and references therein),
whereby the head of the snake can be specified to take any path leaving in its wake a line singularity. 

Before we give a formal argument for this, we introduce a transformed equation and some notation. Writing
$W:=mu^{m-1}$
takes equation \eqref{eq:porous} to the quadratically nonlinear form
\begin{equation} \label{eq:pressure}
W_t=W\Delta W -\frac{1}{1-m}|\nabla W|^2.
\end{equation}
We use the notation $x=(x_1,x_2,\dots,x_n)\in\mathbb{R}^n$, 
$\rho:=(x_2^2+\dots+x_n^2)^{1/2}$.
A formal argument
proceeds along the following lines: at a point on which the singular curve is smooth, we take the $x_1$ direction to be tangential
to the curve, the dominant balance then reading
\begin{equation}\label{eq:bal}
W_t \sim W\left(W_{\rho\rho} +\frac{n-2}{\rho}W_\rho\right)-\frac{1}{1-m}W_\rho^2
\end{equation}
on the assumption of cylindrical symmetry, this having a self-consistent local solution
\begin{equation}\label{eq:locsol}
W \sim \frac{\rho^2}{2B(t-t_*)}\, ,\quad B:=\frac{n-1}{1-m}\left(m-\frac{n-3}{n-1}\right),\qquad\mbox{ as }
\rho\to 0,~ t>t_*\, ,
\end{equation}
where $t_*$ is the time at which the head of the snake passes through the location in question, the behavior at the head 
being a quasi-steady generalization of the above traveling wave solution.

Numerous natural generalizations presumably arise: the head of the snake can come to a halt or retreat, its path may not need to be smooth and
singular sets of dimensionality greater than one are possible (the simplest such examples being the above solutions embedded
in higher dimensional space with no dependence on the additional dimensions). Very specific questions relate to whether cylindrical symmetry
necessarily follows in the sense of \eqref{eq:bal}-\eqref{eq:locsol} and with respect to the large-time behavior when the head comes to a stop in finite time.

\appendix
\section{The traveling wave solution}\label{sec:app}
We seek a traveling wave for \eqref{eq:pressure} in the $x_1$ direction, writing
\[
w=W(x_1-ta, x_2,\dots, x_n),\qquad a>0
\]
and introducing paraboloidal coordinates $Y$ and $Z$ via
\[
x_1-ta=\frac{1}{2}(Y^2-Z^2),\qquad R=YZ, \qquad W=\phi(Y,Z)
\]
to give
\[
\begin{aligned}
-a(Y \phi_Y -Z \phi_Z)&=\phi\left(\phi_{YY} +\phi_{ZZ} +(n-2)\left(\frac{1}{Y} \phi_Y+\frac{1}{Z}\phi_Z\right)\right) \\
&\quad -\frac{1}{1-m}\left(\phi_Y^2+\phi_Z^2\right).
\end{aligned}
\]
Setting $\phi=\phi(Y)$ yields the ODE
\[
-aY \phi_Y =\phi \left(\phi_{YY}  +\frac{n-2}{Y} \phi_Y\right)
-\frac{1}{1-m}\phi_Y^2
\]
that corresponds to the Boltzmann similarity reduction of the porous-medium equation in $n-1$ dimensions and whose
scaling properties imply the existence of a solution $\phi=aY^2/B$ equivalent to the representation \eqref{eq:expli}-\eqref{eq:Adef}.
This solution thus constitutes a reduction 
 of \eqref{eq:porous} akin to that exploited by the Ivantsov solution to the Stefan problem.

\section*{Acknowledgments}
The authors express their gratitude 
to Dr.~Takashi Kagaya for his valuable comments concerning 
the regularity of the distance function and to the referee for useful suggestions.

%\bibliographystyle{amsplain}

\end{document}